\documentclass[a4paper,12pt,twoside]{article}

\usepackage{amssymb}
\usepackage{array}
\usepackage{dsfont}
\usepackage{amsmath}
\usepackage{float}
\usepackage{fancyhdr}
\usepackage[matrix,arrow,curve]{xy}
\usepackage{pstricks} 
\usepackage{amsmath,amsfonts,verbatim,afterpage,theorem,euscript,mathrsfs,amssymb}
\usepackage[english]{babel}
\usepackage{authblk}
 
\newtheorem{theorem}{Theorem}[section]
\newtheorem{lemma}[theorem]{Lemma}
\newtheorem{e-proposition}[theorem]{Proposition}

\newtheorem{e-definition}[theorem]{Definition\rm}

\newtheorem{theoreme}{Th\'eor\`eme}[section]

\newtheorem{proposition}[theoreme]{Proposition}

\def \vg{\overrightarrow{g}}
\def \F{\overrightarrow{F}}
\def \fe{\overrightarrow{f}}
\def \U{\overrightarrow{U}}

\def \vu{\overrightarrow{u}}

\def \P{\mathbb{P}}
\def \Rt{\mathbb{R}^{3}}

\def \finpv{\hfill $\blacksquare$}

\def\PM{\mathcal{PM}}

\title{Frequency decay for Navier-Stokes stationary solutions}

\author[1]{Diego Chamorro}
\author[1]{Oscar Jarr\'in}
\author[1]{Pierre-Gilles Lemari\'e-Rieusset\thanks{Corresponding author: \texttt{plemarie@univ-evry.fr}\\ }}
 
\affil[1]{Laboratoire de Math\'ematiques et Mod\'elisation d'Evry (LaMME), UMR 8071\\ Universit\'e d'Evry Val d'Essonne, 23 Boulevard de France 91037, Evry }

\begin{document}

%\author[authorlabel1]{Diego Chamorro},
%\ead{diego.chamorro@univ-evry.fr}
%\author[authorlabel1]{Oscar Jarr\'in},
%\ead{oscar.jarrin@univ-evry.fr}
%\author[authorlabel1]{Pierre-Gilles Lemari\'e-Rieusset},
%\ead{plemarie@univ-evry.fr}
%\address[authorlabel1]{}

\maketitle

\begin{abstract}
We consider stationary Navier-Stokes equations in $\Rt$ with a regular external force and we prove exponential frequency decay of the solutions. Moreover, if the external force is small enough, we give a pointwise exponential frequency decay for such solutions according to the K41 theory. If a damping term is added to the equation, a pointwise decay is obtained without the smallness condition over the force.
\vskip 0.5\baselineskip
\end{abstract}

%%%%%%%%%%%%%%%%%%%%%%%%%%%%%%%%%%%%%%%%%%%%%%%%%%%%%%%%%%%%%%%%%%%%%%%%%%%%%%%%%%%%%%%%%%%%%%%%%%%%%%%%%%%%%%%%%%%%%

\section{Introduction} 
Gevrey regularity for solutions of the Navier-Stokes equations has been studied in many different frameworks: for a periodic setting with external force see \cite{FoiasTemam}, \cite{Liu}; for the stationary problem in $\mathbb{T}^{3}$ with frequency localized forces see \cite{KLT}. For the evolution problem in $\Rt$ (with a null force) a pointwise analysis is obtained in \cite{PGLR0}.\\

In this article we generalize some of these previous results in the framework of stationary Navier-Stokes equations in $\Rt$
\begin{equation}\label{N-S-stationnaire}
-\nu \Delta \U+\P(div(\U\otimes \U)) =\F, \qquad div(\U)=0,\quad div(\F),
\end{equation}
where $\nu >0$ is the fluid's viscosity parameter, $\U:\mathbb{R}^3\longrightarrow \mathbb{R}^3$ is the velocity, $\P$ is the Leray's projector and $\F:\mathbb{R}^3\longrightarrow \mathbb{R}^3$ is a time-independent external force. \\

If the external force is regular enough we prove in Theorem \ref{Theo_1} an exponential frequency decay. Moreover, if the external force is small enough, we give in Theorem \ref{Theo_3} a pointwise exponential frequency decay for such solutions. Finally, if a damping term is added to the equation, a pointwise decay is obtained in Theorem \ref{Theo_4} without the smallness condition over the force.
%%%%%%%%%%%%%%%%%%%%%%%%%%%%%%%%%%%%%%%%%%%%%%%%%%%%%%%%%%%%%%%%%
\begin{theorem}\label{Theo_1} Let $\F \in \dot{H}^{-1}(\mathbb{R}^3)$ be such that for $\varepsilon_0>0$ we have 
$$\displaystyle{\int_{\mathbb{R}^3}}e^{2\varepsilon_0\vert \xi \vert}\big| \widehat{\F}(\xi)\big|^2 \vert \xi \vert^{-2}d\xi<+\infty.$$
Then there exists $\U\in \dot{H}^{1}(\mathbb{R}^3)$ a solution to the stationary Navier-Stokes equations (\ref{N-S-stationnaire}), such that $\U$ verifies the following exponential frequency decay: 
\begin{equation}\label{frec_decay_U}
\int_{\mathbb{R}^3} e^{2\varepsilon_1\vert \xi \vert}\big|\widehat{\U}(\xi)\big|^2 \vert \xi \vert^2 d\xi<+\infty,\quad \mbox{where $\varepsilon_1=\varepsilon_1(\varepsilon_0,\F,\nu)>0$.}  
\end{equation} 
\end{theorem}
\quad\\[4mm]
%%%%%%%%%%%%%%%%%%%%%%%%%%%%%%%%%%%%%%%%%%%%%%%%%%%%%%%%%%%%%%%%%
In the laminar setting we obtain a sharper \emph{pointwise} exponential frequency decay.\\ For $0\leq a <3$, we define the  pseudo-measures space by 
$$\PM^a=\left\lbrace \vg\in \mathcal{S}^{'}(\mathbb{R}^3): \widehat{\vg}\in L^{1}_{loc}(\mathbb{R}^3)\quad \text{and}\quad \vert\xi \vert^a \widehat{\vg}\in L^{\infty}(\mathbb{R}^3)\right\rbrace,$$ 
which is a Banach space endowed with the norm $\Vert \vg \Vert_{\PM^a}=\Vert \vert \xi \vert^a \widehat{\vg}\Vert_{L^{\infty}}$, for $a=0$ we will simply denote the space $\PM^0$ by $\PM$.  
%%%%%%%%%%%%%%%%%%%%%%%%%%%%%%%%%%%%%%%%%%%%%%%%%%%%%%%%%%%%%%%%%
\begin{theorem}\label{Theo_3} Let  $\F\in \PM$. There exists a (small) constant $\eta>0$ such that if 
$$\underset{\xi \in \mathbb{R}^3}{\sup} e^{\vert \xi \vert}\big|\widehat{\F}(\xi)\big|<\eta,$$ 
then there exists $\U\in \PM^2$ a solution to the stationary Navier-Stokes equations (\ref{N-S-stationnaire}) such that $\U$ verifies the following pointwise exponential frequency decay:
\begin{equation}\label{3hyp_U}
\big|\widehat{\U}(\xi)\big|\leq ce^{-\vert\xi \vert}\vert \xi \vert^{-2},\quad \text{for all}\quad \xi\neq 0.
\end{equation} 
\end{theorem}
%%%%%%%%%%%%%%%%%%%%%%%%%%%%%%%%%%%%%%%%%%%%%%%%%%%%%%%%%%%%%%%%%
If a damping term is added to the stationary Navier-Stokes system, we have the following result

\begin{theorem}\label{Theo_4} Let $\F\in H^{-1}(\Rt)$ and for $\alpha>0$ consider the damped stationary Navier-Stokes equations
\begin{equation}\label{DampedNS}
-\nu \Delta \U+\P(div(\U\otimes \U)) =\F-\alpha\U, \qquad div(\U)=0. 
\end{equation}
If the external force $\F$ is such that $\big|\widehat{\F}(\xi)\big|\leq e^{-\varepsilon_{0}|\xi|}$ for a fixed $\varepsilon_{0}>0$, then the stationary solution $\U\in H^{1}(\Rt)$ satisfies the following pointwise exponential frequency decay
\begin{equation}\label{FreqDecayDamped}
\big|\widehat{\U}(\xi)\big|\leq c e^{-\varepsilon_{1}\vert\xi \vert}\vert \xi \vert^{-\frac{5}{2}},\quad \text{for all}\quad \xi\neq 0, \quad \mbox{where $\varepsilon_1=\varepsilon_1(\varepsilon_0,\F,\nu)>0$.}  
\end{equation} 
\end{theorem}
%%%%%%%%%%%%%%%%%%%%%%%%%%%%%%%%%%%%%%%%%%%%%%%%%%%%%%%%%%%%%%%%%
\section{Proof of Theorem \ref{Theo_1}}
\begin{lemma}\label{Lemma11} If $\F\in \dot{H}^{-1}(\Rt)$, then there exists at least one solution $\U\in \dot{H}^{1}(\Rt)$ to the stationary Navier-Stokes equation (\ref{N-S-stationnaire}).
\end{lemma}
\begin{lemma}\label{Lemma1}
Let $T_0>0$. For $\vu_0 \in \dot{H}^{1}(\mathbb{R}^3)$ a divergence-free initial data and a divergence-free external force $\fe \in \mathcal{C}([0,T_0[,\dot{H}^{1}(\mathbb{R}^3))$   there exists a time $0<T_1<T_0$ and a function $\vu \in \mathcal{C}([0,T_1[,\dot{H}^{1}(\mathbb{R}^3))$ which is a unique solution to the Navier-Stokes equations
\begin{equation}\label{N-S-1}
\partial_t \vu-\nu\Delta \vu+\P(div(\vu\otimes \vu))=\fe,\quad div(\vu)=0,\qquad \vu(0,\cdot)=\vu_0.
\end{equation}  
\end{lemma}
Existence and uniqueness issues are classical, see \cite{PGLR1} for details.\\ 

In the following proposition we prove the frequency decay for the solution $\vu$ obtained in Lemma \ref{Lemma1}.  
%%%%%%%%%%%%%%
\begin{proposition}\label{Lemma2} Let  $\alpha>0$ and consider the Poisson kernel  $e^{\alpha\sqrt{t}\sqrt{-\Delta}}$.  Within the framework  of  Lemma \ref{Lemma1}, if the external force   $\fe$ is  such that 
$$e^{\alpha\sqrt{t}\sqrt{-\Delta}} \fe \in \mathcal{C}(]0,T_0[,\dot{H}^{1}(\mathbb{R}^3)),$$ 
then the unique solution of  equations (\ref{N-S-1}) satisfies $e^{\alpha\sqrt{t}\sqrt{-\Delta}} \vu \in \mathcal{C}(]0,T_1[,\dot{H}^{1}(\mathbb{R}^3))$ for all time $t\in [0,T_1[$ where $0<T_1<T_0$ is small enough.  
\end{proposition}
%%%%%%%%%%%%%%
\textbf{Proof.}  Consider the space 
$$E=\left\lbrace \vu \in \mathcal{C}(]0,T_1[,\dot{H}^{1}(\mathbb{R}^3)): e^{\alpha\sqrt{t}\sqrt{-\Delta}} \vu \in \mathcal{C}(]0,T_1[,\dot{H}^{1}(\mathbb{R}^3)) \right\rbrace,$$ 
endowed with the norm $\Vert \cdot \Vert_E= \Vert e^{\alpha\sqrt{t}\sqrt{-\Delta}}(\cdot)\Vert_{L^{\infty}_{t}\dot{H}^{1}_{x}}$. We study the quantity   
\begin{equation}\label{EstimationPointFixe1}
\Vert \vu_1 \Vert_{E}= \left\Vert h_{\nu t}\ast \vu_0+\int_{0}^{t}h_{\nu(t-s)}\ast \fe(s,\cdot)ds-\int_{0}^{t}h_{\nu(t-s)}\ast \P(div(\vu_1 \otimes \vu_1))(s,\cdot)ds\right\Vert_{E}
\end{equation} 
where $h_{\nu t}$ is the heat kernel. The two first terms of this expression are easy to estimate and we have
\begin{equation}\label{eq04}
\left\Vert h_{\nu t}\ast \vu_0+\int_{0}^{t}h_{\nu(t-s)}\ast \fe(s,\cdot)ds\right\Vert_{E}\leq c(\nu,\alpha,T_0)\left(\Vert \vu_0\Vert_{\dot{H}^{1}_{x}}+\Vert e^{\alpha\sqrt{t}\sqrt{-\Delta}}\fe \Vert_{L^{\infty}_{t}\dot{H}^{1}_{x}}\right).
\end{equation}
For the last term of (\ref{EstimationPointFixe1}), by definition of the norm $\Vert \cdot \Vert_E$, by the Plancherel formula and by the boundedness of the Leray projector we have 
\begin{eqnarray*}
(I)&=&\left\Vert \int_{0}^{t}h_{\nu(t-s)}\ast \P(div(\vu_1 \otimes \vu_1))ds\right\Vert_{E}\\
& = &\sup_{0<t<T_1} \left\Vert e^{\alpha\sqrt{t}\sqrt{-\Delta}}\left(\int_{0}^{t}h_{\nu(t-s)}\ast \P(div(\vu_1 \otimes \vu_1)) ds \right)\right\Vert_{\dot{H}^{1}_{x}}\\
&\leq &\sup_{0<t<T_1} c \left\Vert \vert\xi \vert^2 \int_{0}^{t}e^{-\nu(t-s)\vert \xi \vert^2} e^{\alpha\sqrt{t}\vert\xi \vert}\left\vert\left(\mathcal{F}[\vu_1]\ast \mathcal{F}[\vu_1]\right)(s,\cdot) \right\vert ds\right\Vert_{L^{2}_{x}}.
\end{eqnarray*}
Since we have the pointwise inequality 
\begin{equation}\label{eq10}
e^{\alpha\sqrt{t}\vert\xi \vert}\left\vert\left(\mathcal{F}[\vu_1]\ast \mathcal{F}[\vu_1]\right)(s,\xi)\right\vert \leq \left[\left( e^{\alpha\sqrt{t}\vert\xi \vert}\vert \mathcal{F}[\vu_1]\vert\right)\ast \left( e^{\alpha\sqrt{t}\vert\xi \vert}\vert \mathcal{F}[\vu_1]\vert\right)\right](s,\xi),
\end{equation}
due to the fact that $e^{\alpha\sqrt{t}\vert \xi \vert}\leq e^{\alpha\sqrt{t}\vert \vert\xi-\eta\vert}e^{\alpha\sqrt{t}\vert \eta \vert}$ for all $\xi,\,\eta\in \mathbb{R}^3$, 
then we obtain
\begin{eqnarray*}
(I)&\leq & \sup_{0<t<T_1} c \int_{0}^{t}\left\Vert \vert\xi \vert^\frac{3}{2} e^{-\nu(t-s)\vert \xi \vert^2} \vert\xi \vert^{\frac{1}{2}}\left\vert\left[\left( e^{\alpha\sqrt{t}\vert\xi \vert}\vert \mathcal{F}[\vu_1]\vert\right)\ast \left( e^{\alpha\sqrt{t}\vert\xi \vert}\vert \mathcal{F}[\vu_1]\vert\right)\right]\right\vert \right\Vert_{L^{2}_{x}}ds.
\end{eqnarray*}
Getting back to the spatial variable we can write
\begin{eqnarray}
(I)&\leq &\sup_{0<t<T_1} c \int_{0}^{t} \left\Vert (-\Delta)^{\frac{3}{4}}h_{\nu(t-s)}\ast  (-\Delta)^{\frac{1}{4}} \left\lbrace\left(\mathcal{F}^{-1}\left[e^{\alpha\sqrt{t}\vert\xi \vert}\vert \mathcal{F}[\vu_1]\vert \right]\right) \otimes \right. \right.\nonumber\\
& &\qquad \qquad \left.\left. \left(\mathcal{F}^{-1}\left[e^{\alpha\sqrt{t}\vert\xi \vert}\vert \mathcal{F}[\vu_1]\vert\right]\right)\right\rbrace\right\Vert_{L^{2}_{x}}ds\nonumber\\
&\leq & \left(c  \int_{0}^{T_1} \left\Vert (-\Delta)^{\frac{3}{4}}h_{\nu(t-s)} \right\Vert_{L^1}ds\right) \left\Vert  \left(\mathcal{F}^{-1}\left[e^{\alpha\sqrt{t}\vert\xi \vert}\vert \mathcal{F}[\vu_1]\vert \right]\right) \otimes \right.\nonumber\\
& &\qquad \qquad   \left. \left(\mathcal{F}^{-1}\left[e^{\alpha\sqrt{t}\vert\xi \vert}\vert \mathcal{F}[\vu_1]\vert\right]\right) \right\Vert_{L^{\infty}_{t}\dot{H}^{\frac{1}{2}}_{x}}\nonumber\\
&\leq& c\frac{T^{\frac{1}{4}}}{\nu^{\frac{3}{4}}} \left\Vert \mathcal{F}^{-1}\left[e^{\alpha\sqrt{t}\vert\xi \vert}\vert \mathcal{F}[\vu_1]\vert \right]\right\Vert_{L^{\infty}_{t}\dot{H}^{1}_{x}} \left\Vert \mathcal{F}^{-1}\left[e^{\alpha\sqrt{t}\vert\xi \vert}\vert \mathcal{F}[\vu_1]\vert \right]\right\Vert_{L^{\infty}_{t}\dot{H}^{1}_{x}}\nonumber\\
&\leq &c\frac{T^{\frac{1}{4}}_{1}}{\nu^{\frac{3}{4}}}\Vert \vu_1 \Vert_E \Vert \vu_1 \Vert_E.\label{eq08}
\end{eqnarray} 
With estimates (\ref{eq04}) and (\ref{eq08}) at hand, we fix $T_1$ small enough in order to apply Picard's contraction principle and we obtain a solution $\vu_1\in E$  of (\ref{N-S-1}). Since $E\subset \mathcal{C}(]0,T_1[,\dot{H}^{1}(\mathbb{R}^3))$ we have $\vu_1\in \mathcal{C}(]0,T_1[,\dot{H}^{1}(\mathbb{R}^3))$ and by uniqueness of the solution $\vu$ we have $\vu_1=\vu$, and thus $\vu \in E$. \finpv\\
\\
Now, we come back to the stationary Navier-Stokes equations (\ref{N-S-stationnaire}) and we will prove that the solution $\U\in \dot{H}^{1}(\mathbb{R}^3)$ (given by Lemma \ref{Lemma11}) satisfies the exponential frequency decay given in (\ref{frec_decay_U}). In the space $\mathcal{C}(]0,1[,\dot{H}^{1}(\mathbb{R}^3))$ we consider the evolution problem (\ref{N-S-1}) with the initial data $\vu_{0}=\U$ where the external force $\fe$ is now given by with the expression 
$$ \fe = e^{-\alpha\sqrt{t}\sqrt{-\Delta}}(e^{\alpha \sqrt{t}\sqrt{-\Delta}}\F),$$ 
for the particular value $\alpha =\frac{2}{3}\varepsilon_0>0$ where $\varepsilon_0>0$ is given in the hypothesis of the force $\F$. To obtain a unique solution $\vu\in \mathcal{C}(]0,1[,\dot{H}^{1}(\mathbb{R}^3))$ to the equations (\ref{N-S-1}) such that 
$$e^{\alpha\sqrt{t}\sqrt{-\Delta}}\vu\in \mathcal{C}(]0,1[,\dot{H}^{1}(\mathbb{R}^3)),$$ 
we prove that the external force  $\fe$ verifies the hypotheses of Lemma \ref{Lemma1} and Proposition \ref{Lemma2} above: 
\begin{eqnarray*}
\left\Vert e^{\alpha\sqrt{t}\sqrt{-\Delta}}\F \right\Vert^{2}_{L^{\infty}_{t}\dot{H}^{1}_{x}} &=&\sup_{0<t<1} \int_{\mathbb{R}^3} \vert \xi \vert^2e^{2\alpha\sqrt{t} \vert \xi \vert}\big|\widehat{\F}(\xi)\big|^2d\xi \\
&\leq & \frac{1}{\alpha^4}\int_{\mathbb{R}^3} (\alpha\vert \xi \vert)^4 e^{2\alpha \vert \xi \vert}\big| \widehat{\F}(\xi)\big|^2\, \vert \xi \vert^{-2}d\xi\\
&\leq & \frac{1}{\alpha^4} \int_{\mathbb{R}^3} e^{3\alpha \vert \xi \vert}\big| \widehat{\F}(\xi)\big|^2\, \vert \xi \vert^{-2}d\xi\\
&\leq &  \frac{1}{\alpha^4} \int_{\mathbb{R}^3} e^{2\varepsilon_0\vert \xi \vert}\big|\widehat{\F}(\xi)\big|^2\, \vert \xi \vert^{-2}d\xi<+\infty.\\
\end{eqnarray*} Thus, once we have $e^{\alpha\sqrt{t}\sqrt{-\Delta}}\F\in \mathcal{C}(]0,1[,\dot{H}^{1}(\mathbb{R}^3))$, since the operator $e^{-\alpha\sqrt{t}\sqrt{-\Delta}}$ is bounded in the space $\mathcal{C}(]0,1[,\dot{H}^{1}(\mathbb{R}^3))$ we have  
$$\fe=e^{-\alpha\sqrt{t}\sqrt{-\Delta}}(e^{\alpha\sqrt{t}\sqrt{-\Delta}}\F)\in \mathcal{C}(]0,1[,\dot{H}^{1}(\mathbb{R}^3)).$$
Moreover, we have 
$$ e^{\alpha\sqrt{t}\sqrt{-\Delta}}\fe =e^{\alpha\sqrt{t}\sqrt{-\Delta}}\F\in \mathcal{C}(]0,1[,\dot{H}^{1}(\mathbb{R}^3)).$$ 
By  Lemma \ref{Lemma1}  there exists a time $0<T_1<1$ and a unique solution $\vu \in \mathcal{C}(]0,T_1[,\dot{H}^{1}(\mathbb{R}^3))$ to the equation (\ref{N-S-1}). Moreover, since $e^{\alpha\sqrt{t}\sqrt{-\Delta}}\fe\in \mathcal{C}(]0,1[,\dot{H}^{1}(\mathbb{R}^3))$ by Proposition \ref{Lemma2} we have  $ e^{\alpha\sqrt{t}\sqrt{-\Delta}}\vu \in \mathcal{C}(]0,T_1[,\dot{H}^{1}(\mathbb{R}^3))$. Since the solution $\U\in \dot{H}^{1}(\mathbb{R}^3)$ of the stationary Navier-Stokes equations (\ref{N-S-stationnaire}) is a constant in time, we have  $\U\in \mathcal{C}(]0,T_1[,\dot{H}^{1}(\mathbb{R}^3))$ and since $\partial_t\U\equiv 0$ and 
 $$\fe =e^{-\alpha\sqrt{t}\sqrt{-\Delta}}(e^{\alpha\sqrt{t}\sqrt{-\Delta}}\F) =\F,$$ 
 we find that $\U \in \mathcal{C}(]0,T_1[,\dot{H}^{1}(\mathbb{R}^3))$ is also a solution to the equation (\ref{N-S-1}) and thus,  by uniqueness we get  $\U=\vu$. Then, since  $e^{\alpha\sqrt{t}\sqrt{-\Delta}}\vu \in \mathcal{C}(]0,T_1[,\dot{H}^{1}(\mathbb{R}^3))$ we have  
 $$e^{\alpha\sqrt{t}\sqrt{-\Delta}}\U \in \mathcal{C}(]0,T_1[,\dot{H}^{1}(\mathbb{R}^3)),$$ 
 for all time $t\in [0,T_1[$. Thus, if $\varepsilon_1=\alpha\sqrt{\frac{T_1}{2}} >0$, we have 
\begin{equation*}
 \int_{\mathbb{R}^3}e^{2\varepsilon_1\vert \xi \vert} \vert \U(\xi)\vert^2 \vert \xi \vert^2 d\xi =\big\| e^{\alpha\sqrt{\frac{T_1}{2}}\sqrt{-\Delta}}\U\big\|^{2}_{\dot{H}^{1}_{x}}\leq \sup_{0<t<T_1}\Vert e^{\alpha\sqrt{t}\sqrt{-\Delta}}\U\Vert^{2}_{\dot{H}^{1}_{x}}<+\infty,
\end{equation*} 
and we obtain the frequency decay given in (\ref{frec_decay_U}). \finpv
%%%%%%%%%%%%%%%%%%%%%%%%%%%%%%%%%%%%%%%%%%%%%%%%%%%%%%%%%%%%%%%%%%%%%%%%%%%%%%%%%%%%%%%%%%%%%%%%%%%%%%%%%%%%%%%%%%%%%%%%%%%%%%%%%%%%%%%%%%%%%%%%%%%%%%%%%
\section{Proof of Theorem \ref{Theo_3}}
We consider now the space $A= \left\lbrace \U \in \PM^2: e^{\sqrt{-\Delta}}\U\in \PM^2\right\rbrace$, endowed with the norm
\begin{equation}\label{3normE}
\Vert \cdot\Vert_A=\Vert e^{\sqrt{-\Delta}}(\cdot)\Vert_{\PM^2},
\end{equation}
and in this space we study the existence of a solution of equations (\ref{N-S-stationnaire}) under the hypotheses of Theorem \ref{Theo_3}. For this we study the quantity
\begin{equation}\label{3eq01}
\big\| \U\big\|_{A}=\left\Vert \frac{1}{\nu}\P\left(\frac{1}{\Delta}div(\U\otimes \U)\right)-\frac{1}{\nu}\frac{1}{\Delta}\F\right\Vert_{A} \leq  \frac{1}{\nu}\left\Vert \P\left(\frac{1}{\Delta}div(\U\otimes \U)\right)\right\Vert_{A}+\frac{1}{\nu}\left\Vert \frac{1}{\Delta}\F\right\Vert_{A},
\end{equation} where, for the first term of the inequality above we have the following estimate:
\begin{equation}\label{3estimate_B}
\frac{1}{\nu}\left\Vert\P\left(\frac{1}{\Delta}div(\U\otimes \U)\right)\right\Vert_{A}\leq \frac{c}{\nu}\Vert \U\Vert_{A}\Vert \U\Vert_{A}. 
\end{equation}
Indeed, by the expression (\ref{3normE}) and by the continuity of the Leray projector we have 
\begin{eqnarray}
\frac{1}{\nu}\left\Vert \P\left(\frac{1}{\Delta}div(\U\otimes \U)\right)\right\Vert_{A} &=&
\frac{1}{\nu}\left\Vert \vert\xi \vert^{2} e^{\vert \xi \vert} \mathcal{F}\left[ \P\left(\frac{1}{\Delta}div(\U\otimes \U)\right) \right]\right\Vert_{L^{\infty}}\nonumber\\
& \leq & \frac{c}{\nu} \left\Vert \vert\xi \vert^{2} e^{\vert \xi \vert} \frac{1}{\vert \xi \vert} \left\vert\mathcal{F}\left[ \U\right] \ast\mathcal{F}\left[ \U\right] \right\vert \right\Vert_{L^{\infty}}\nonumber\\
&\leq &\frac{c}{\nu} \left\Vert \vert\xi \vert \left[ \left( e^{\vert \xi \vert}\mathcal{F}\left[ \vert \U \vert \right]\right)\ast \left( e^{\vert \xi \vert}\mathcal{F}\left[ \vert \U \vert \right]\right)  \right]\right\Vert_{L^{\infty}},\label{3estimateB_1}
\end{eqnarray}
where the last inequality can be deduced from (\ref{eq10}). Now we remark that 
\begin{eqnarray*}
\left[ \left( e^{\vert \xi \vert}\mathcal{F}\left[ \vert \U \vert \right]\right)\ast \left( e^{\vert \xi \vert}\mathcal{F}\left[ \vert \U \vert \right]\right)  \right](\xi)&=&\int_{\mathbb{R}^3} e^{\vert \xi-\eta \vert}\mathcal{F}\left[ \vert \U \vert \right](\xi-\eta)e^{\vert \eta \vert}\mathcal{F}\left[ \vert \U \vert \right](\eta)d \eta\\  
&\leq &  \Vert \U \Vert_A\,\Vert \U \Vert_A \int_{\mathbb{R}^3} \frac{d \eta}{\vert \xi-\eta\vert^2 \vert \eta \vert^2}\leq \frac{c}{\vert \xi \vert} \Vert \U \Vert_A\,\Vert \U \Vert_A,
\end{eqnarray*}
and thus, using this inequality in (\ref{3estimateB_1}) we easily obtain the estimate  (\ref{3estimate_B}). For the second term in the RHS of (\ref{3eq01}) we have 
$$\frac{1}{\nu}\left\Vert \frac{1}{\Delta}\F\right\Vert_A = \frac{1}{\nu}\left\Vert e^{\sqrt{-\Delta}}\left(\frac{1}{\Delta}\F \right)\right\Vert_{\PM^2}= \frac{c_1}{\nu}\sup_{\xi \in\mathbb{R}^3} \vert \xi \vert^2 e^{\vert \xi \vert} \frac{1}{\vert \xi \vert^2}\vert \F(\xi)\vert=\frac{c_1}{\nu} \sup_{\xi \in\mathbb{R}^3} e^{\vert \xi \vert} \vert \F(\xi)\vert.$$ 
Thus, if the external force $\F$ satisfies $ \underset{\xi \in\mathbb{R}^3}{\sup} e^{\vert \xi \vert} \vert \F(\xi)\vert<\eta$, for $\eta$ small enough, we obtain $\U\in A$ a solution to the stationary Navier-Stokes equations (\ref{N-S-stationnaire}) for which we have the pointwise estimate (\ref{3hyp_U}). 
\finpv 
%%%%%%%%%%%%%%%%%%%%%%%%%%%%%%%%%%%%%%%%%%%%%%%%%%%%%%%%%%%%%%%%%%%%%%%%%%%%%%%%%%%%%%%%%%%%%%%%%%%%%%%%%%%%%%%%%%%%%%%%%%%%%%%%%%%%%%%%%%%%%%%%%%%%%%%%%
\section{Proof of Theorem \ref{Theo_4}}
For $\alpha>0$ and under the hypotheses of Theorem \ref{Theo_4}, the existence of solutions of equation (\ref{DampedNS}) is given by applying the Scheafer fixed point theorem. Now, for $\vu_0\in \PM^{\frac{5}{2}}$  we consider the non-stationary damped Navier-Stokes equations 
\begin{equation}\label{N-S}
\partial_t \vu+\P(div(\vu\otimes \vu))-\nu\Delta \vu=\fe -\alpha\vu, \quad div(\vu)=0,\quad \vu(0,\cdot)=\vu_0,
\end{equation}
where the divergence-free external force $\fe$ belongs to the space  $\mathcal{C}([0,T_0[,\PM^{\frac{5}{2}})$. For this problem there exists a unique solution $\vu \in \mathcal{C}([0,T_1[,\PM^{\frac{5}{2}})$ with $0<T_{1}<T_{0}$. For existence issues for equations (\ref{DampedNS}) and (\ref{N-S}) see the details in \cite{PGLR1}. \\

Following essentially the same lines of Proposition \ref{Lemma2} above, we prove that if the external force is such that $e^{\beta\sqrt{t}\sqrt{-\Delta}}\fe \in \mathcal{C}([0,T_0[,\PM^{\frac{5}{2}})$  then the unique solution of (\ref{N-S}) is such that $e^{\beta\sqrt{t}\sqrt{-\Delta}}\vu \in \mathcal{C}([0,T_1[,\PM^\frac{5}{2})$. As in the proof of Theorem \ref{Theo_3}, we consider 
$$\fe = e^{-\beta\sqrt{t}\sqrt{-\Delta}}\left(e^{\beta\sqrt{t}\sqrt{-\Delta}}\F\right)=\F,$$ 
and for a suitable value of the parameter $\beta>0$ we can prove that $\fe \in \mathcal{C}([0,1[,\PM^{\frac{5}{2}})$ and $e^{\beta\sqrt{t}\sqrt{-\Delta}}\fe \in \mathcal{C}([0,1[,\PM^{\frac{5}{2}})$.\\ 

In order to link the stationary solution to the non-stationary problem, we must prove that the solution $\U\in H^{1}(\Rt)$ of (\ref{DampedNS}) is such that $\U\in \PM^{\frac{5}{2}}$, and in this step we use the extra damping term. Indeed, rewriting (\ref{DampedNS}) we consider the equation
\begin{equation}\label{int_N-S-stat-damped}
\U=\frac{-\nu \Delta}{\alpha Id -\nu \Delta} \left( \P\left( \frac{1}{\nu \Delta}div(\U\otimes \U) \right) \right)+\frac{1}{\alpha Id -\nu \Delta}\left( \F\right),
\end{equation}
and we obtain 
\begin{eqnarray*}
\Vert \U \Vert_{\dot{H}^{\frac{3}{2}}}\leq & \left\Vert \frac{-\nu \Delta}{\alpha I_d -\nu \Delta} \left( \P\left( \frac{1}{\nu \Delta}div(\U\otimes \U) \right) \right) \right\Vert_{\dot{H}^{\frac{3}{2}}}+ \left\Vert \frac{1}{\alpha I_d -\nu \Delta} \left( \F \right) \right\Vert_{\dot{H}^{\frac{3}{2}}}.\\
\end{eqnarray*}
Since the operator $\frac{-\nu \Delta}{\alpha I_d -\nu \Delta}$ is bounded in $\dot{H}^{\frac{3}{2}}(\Rt)$ and by the properties of $\F$ we can write
\begin{eqnarray*}
\Vert \U \Vert_{\dot{H}^{\frac{3}{2}}}& \leq & \left\Vert \frac{1}{\nu\Delta}div (\U \otimes \U) \right\Vert_{\dot{H}^{\frac{3}{2}}}+\left\Vert  \frac{1}{\alpha I_d -\nu \Delta} \left( \F \right) \right\Vert_{H^{2}}\\ 
&\leq & c\| \U \otimes \U \|_{\dot{H}^{\frac{1}{2}}}+ c(\alpha) \|  \F  \|_{L^2}\\
&\leq & c \Vert \U \Vert_{H^1} \Vert \U \Vert_{H^1}+ c(\alpha) \| \F \|_{L^2}.
\end{eqnarray*}
We thus have $\U\in \dot{H}^{\frac{3}{2}}(\Rt)$ and we prove now $\U \in \PM^{\frac{5}{2}}$:  from equation (\ref{int_N-S-stat-damped}) we obtain 
\begin{eqnarray*}
\left\vert \mathcal{F}\left[\U\right](\xi)\right\vert &\leq & c\frac{1}{\nu \vert \xi \vert} \vert \left(\mathcal{F}\left[ \U \right] \ast \mathcal{F}\left[ \U \right]\right)(\xi)\vert + \frac{1}{\nu \vert \xi \vert^2}\vert \mathcal{F}\left[ \F \right](\xi)\vert,
\end{eqnarray*} and then, multiplying by $\vert \xi \vert^{\frac{5}{2}}$ and by hypothesis on $\F$ we get the estimate
\begin{eqnarray*}
\vert \xi \vert^{\frac{5}{2}} \left\vert \mathcal{F}\left[\U\right](\xi)\right\vert & \leq &  \int_{\Rt}  \vert \xi \vert^{\frac{3}{2}}\left\vert\mathcal{F}\left[\U\right](\xi-\eta)\right\vert \left\vert\mathcal{F}\left[\U\right](\eta)\right\vert d\eta + \frac{1}{\nu}\vert \xi \vert^{\frac{1}{2}}\left\vert \mathcal{F}\left[ \F \right](\xi)\right\vert\\
&\leq & 2 \Vert \U \Vert_{\dot{H}^{\frac{3}{2}}}\Vert \U \Vert_{L^2}+ \frac{1}{\nu} \vert \xi \vert^{\frac{1}{2}} e^{-\varepsilon_0 \vert \xi \vert},
\end{eqnarray*}
from which we deduce that $\U\in \PM^{\frac{5}{2}}$. Then, we study (\ref{N-S}) with $\vu_{0}=\U$ and we have $\U\in \mathcal{C}([0,T_1[,\PM^{\frac{5}{2}})$, but since $\U$ verifies the equations (\ref{DampedNS}), $\partial_t \U\equiv 0$ and $\fe =\F$, we obtain that $\U$ is also a  solution of (\ref{N-S}) and by uniqueness we have $\U=\vu$. Finally, we have $e^{\beta\sqrt{t}\sqrt{-\Delta}}\U \in \mathcal{C}([0,T_1[,\PM^{\frac{5}{2}})$ for $0<t<T_1$ and if $\varepsilon_1=\beta \sqrt{\frac{T_1}{2}}$ we can write 
$$ \|e^{\varepsilon_1\sqrt{-\Delta}} \U \|_{\PM^{\frac{5}{2}}} \leq \|e^{\beta\sqrt{t}\sqrt{-\Delta}}\U\|_{L^{\infty}([0,T_1[,\PM^{\frac{5}{2}})}<+\infty,$$ 
and we obtain the frequency decay stated in the formula (\ref{FreqDecayDamped}). \finpv

\end{document}